\documentclass[10pt]{amsart}
\newcommand{\bracenom}{\genfrac{\lbrace}{\rbrace}{0pt}{}}
\begin{document}
%\doi{10.1080/10652460YYxxxxxxx}
 %\issn{1476-8291}
%\issnp{1065-2469}
% \jvol{00} \jnum{00} \jyear{2024} \jmonth{June}

\markboth{S. Yakubovich}{Integral Transforms and Special Functions}

%\articletype{GUIDE}

 \begin{center} {\huge Uniform upper bounds and asymptotic expansions for certain index kernels}\end{center}
\vspace{2mm}

\begin{center}{{\bf S.Yakubovich}}\end{center}

%\\ \thanks{Email: syakubov@fc.up.pt}

\vspace{6pt}  {\rm {Department of Mathematics, Faculty of Sciences,  University of Porto,  Campo Alegre str.,  687; 4169-007 Porto,  Portugal}\\ \vspace{6pt}}

%\maketitl

\begin{abstract} {\noindent We continue to establish uniform upper bounds and asymptotic expansions  for the kernels of the index transforms which were recently developed for the Kontorovich-Lebedev operator. It involves the Mehler-Fock, Lebedev,  index Whittaker and Olevskii transforms.  }

\end{abstract}
\vspace{4mm}



\vspace{4mm}

\section{Upper bounds for index kernels}

In the theory of the index transforms [1], where the indices or parameters of the hypergeometric functions are happen to be vital components of  the kernels, it is  crucial   to  study  asymptotic behavior at infinity by the index for bounded arguments   or to find upper bounds for their kernels. Recently [2], the author investigated the kernel of the Kontorovich-Lebedev transform,  being the modified Bessel function $K_{i\tau}(x)$ of the pure imaginary index $i\tau$ and positive argument $ x$. It was  proved  the following uniform bound  

$$ \left| K_{i\tau}(x) \right| \le   { \Gamma\left(2^{-n-1}\right)\over  2^{1- 2^{-n}} } \  \bigg[  \sqrt  x \sinh \left(2^{n} \pi\tau /2\right)\bigg]^{- 2^{-n}},  \quad n \in \mathbb{N},\quad x, \tau > 0,\eqno(1.1)$$
where $\Gamma(z)$ is the Euler gamma function. It contains as a particular case the Lebedev inequality  (cf. [1], p. 219).  A precise asymptotic expansion by index $\tau \to +\infty$ and a fixed positive argument is given by 

{\bf Theorem 1\ (cf. [2])}. {\it Let $N$ be a non-negative  integer and $ x \in \left(0, X\right],\ X >0$. Then the modified Bessel function $K_{i\tau}(x)$  has the following asymptotic expansion

$$ K_{i\tau}(x)=  \sqrt{{2\pi\over \tau}}\  e^{-\pi\tau/2} \bigg[ \cos \left( \tau \log \left({2\tau\over e x}\right)  - {\pi\over 4} \right) + R_N (\tau) \bigg],\quad \tau \to +\infty,\eqno(1.2)$$
where the remainder term is given explicitly

$$ R_N (\tau) =  {\rm Re} \bigg[ \exp \left( i \left( \tau \log \left({2\tau\over e x}\right) - {\pi\over 4} \right) \right) \bigg[ r(i\tau) + \bigg( 1+ r(i\tau) \bigg)  \bigg.\bigg.$$

$$ \times \left[ \sum_{m=1}^{N} {(x/2)^{2m}\over  m! \ (1-i\tau)_m} + { x^{2(N+1)} \over 4^{N+1} (N+1)!  (1-i\tau)_{N+1}} \right. $$

$$\left.\left.\left.  \times  {}_1F_2 \left( 1;\ N-i\tau+2,\ N +2;\ {x^2\over 4}\right)\right] \right] \right],\eqno(1.3)$$
and 

$$r(i\tau) = \exp\left( \int_0^\infty e^{-i\tau t} \left[ {1\over 2}- {1\over t} + {1\over e^t-1}\right] {dt\over t} \right) -1.\eqno(1.4)$$
Moreover, the remainder term $R_N(\tau)$ has the following upper bound

$$  \left| R_N (\tau) \right| \le  {1\over \tau} \left[ {e^{1/(6\tau_0)} \over 6}+  \left( \tau_0 +  { e^{1/(6\tau_0)} \over 6 }\right)  \left[ \exp \left( {X^2\over 4\tau_0} \right) \right.\right.$$

$$\left.\left.+ \left( {X^2 \over 2\tau_0}\right)^N  \left[ {I_{N}\left( X \right) \over X^N} - {1\over 2^N N!} \right]\right] \right],\eqno(1.5)$$
where $\tau \ge \tau_0 > 0$ and $I_\nu(x)$ is the modified Bessel function of the first kind.}

{\bf Remark 1}. Comparing  (1.3) with (2.9) in [2] for the remainder term $R_N(\tau)$, the corresponding integral is calculated via Entry 2.15.2.5 in [3], Vol. II in terms of the hypergeometric function ${}_1F_2$ which yields the value

$$  \int_0^x (x^2-y^2)^{N-i\tau} I_{N+1}\left(  y \right) {dy\over y^N} = {2^{-N-2} x^{2N -2i\tau +2} \over (N-i\tau+1) (N+1)!}$$

$$\times  {}_1F_2 \left( 1;\ N-i\tau+2,\ N +2;\ {x^2\over 4}\right).$$

Our goal here is to establish similar results for  the Mehler-Fock, index Whittaker, Olevskii transforms and Lebedev's transforms with the square and the product of the modified Bessel functions.   Indeed,  to estimate the kernel of the Mehler-Fock transform [1], which contains the associated Legendre function $P_\nu^\mu(z)$ [3], Vol. III, we employ the following integral representation (see [3], Vol. II, Entry 2.16.37.3)

$$ \left|\Gamma\left(\mu+ {1\over 2}+i\tau\right)\right|^2  \left[ P^{-\mu}_{-1/2+i\tau} \left(\sqrt{1+4x^2}\right)\right]^2$$

$$ = 2 \int_0^\infty \left[J_\mu(xy)\right]^2 K_{2i\tau}\left(y\right) dy,\quad  \mu > - {1\over 2}. \eqno(1.6)$$
Hence, recalling (1.1)  and Entry 2.12.31.2 in [3], Vol. II, we find 

$$\left| P^{-\mu}_{-1/2+i\tau} \left(\sqrt{1+4x^2}\right)\right| \le  2^{ (2^{-n}+ 1)/2} \  \Gamma^{1/2} \left(2^{-n-1}\right)\  $$

$$ \times \   x^{(2^{-n} -1)/2} \left(\int_0^\infty \left[J_\mu( y)\right]^2 y^{- 2^{-n}} dy\right)^{1/2}  { \bigg[ \sinh \left(2^{n} \pi\tau \right)\bigg]^{- 2^{-n-1}} \over \left|\Gamma\left(\mu+ 1/2+i\tau\right)\right|}$$

$$=   {2^{2^{-n-1}}\over \pi^{1/4}}   \Gamma \left(2^{-n-1}\right)\   \bigg[ { \Gamma (1/2+ \mu -2^{-n-1}) \over  \Gamma ((1+2^{-n })/2) \Gamma (1/2+ \mu+ 2^{-n-1 }) } \bigg]^{1/2} $$

$$ \times  { \bigg[ \sinh \left(2^{n} \pi\tau \right)\bigg]^{- 2^{-n-1}}   x^{(2^{-n} -1)/2} \over \left|\Gamma\left(\mu+ 1/2+i\tau\right)\right|}.$$
Thus we establish the following inequality

$$\left| P^{-\mu}_{-1/2+i\tau} \left(\sqrt{1+4x^2}\right)\right| \le  {2^{2^{-n-1}}\over \pi^{1/4}}   \Gamma \left(2^{-n-1}\right)\   \bigg[ { \Gamma (1/2+ \mu -2^{-n-1}) \over  \Gamma (1/2 +2^{-n-1 }) \Gamma (1/2+ \mu+ 2^{-n-1 }) } \bigg]^{1/2} $$

$$ \times \  { \bigg[ \sinh \left(2^{n} \pi\tau \right)\bigg]^{- 2^{-n-1}}   x^{(2^{-n} -1)/2} \over \left|\Gamma\left(\mu+ 1/2+i\tau\right)\right|},\quad x, \tau > 0,\ n \in \mathbb{N}.\eqno(1.7)$$
Further, the Lebedev index transform with the product of the modified Bessel functions is related with the kernel $[I_{i\tau}(x)+  I_{-i\tau}(x)] K_{i\tau} (x)$.  It has the integral representation (see [3], Vol. II, Entry 2.12.26.1)

$$[I_{i\tau}(x)+  I_{-i\tau}(x)] K_{i\tau} (x) = 2  \int_0^\infty J_0(2x \sinh t) \cos(2 \tau t) dt.\eqno(1.8)$$
Hence we easily find

$$\left| [I_{i\tau}(x)+  I_{-i\tau}(x)] K_{i\tau} (x) \right| \le  2\sqrt 2   \sup_{t\ge 0} \sqrt t\  |J_0(t)| \ x^{-1/2} \int_0^\infty {dy\over  (y^4+1)^{1/2} } $$

$$= {1\over \sqrt {2 \pi}}  \sup_{t\ge 0} \sqrt t\  |J_0(t)| \   \Gamma^2 \left({1\over 4}\right) \ x^{-1/2} .$$
Taking into account the value of the letter supremum (cf. [2]), it gives the bound

$$ \left| [I_{i\tau}(x)+  I_{-i\tau}(x)] K_{i\tau} (x) \right|  \le    \pi^{-1} \  \Gamma^2 \left({1\over 4}\right)\ x^{-1/2},\quad x, \tau > 0.\eqno(1.9)$$
The index transform with the Whittaker function $W_{\mu, i\tau}(x)$ [3], Vol. III as the kernel (index Whittaker transform) can be treated via the representation (cf. [3], Vol. II, Entry 2.16.7.3) 

$$W_{- \mu,i\tau} (2x)=  {1\over \Gamma(\mu)} \sqrt {{2x\over \pi}}  \int_0^\infty  {y^{\mu-1} e^{-  x y} \over (y+1)^{\mu+1/2}}\    K_{i\tau} (x(y+1)) dy,\eqno(1.10)$$
where $x,\ \tau  > 0, \ {\rm Re} \mu > 0.$  Hence,  employing again the generalized Lebedev inequality (1.1), we obtain

$$\left| W_{- \mu,i\tau} (2x) \right| \le  {\Gamma\left(2^{-n-1}\right)  \Gamma\left({\rm Re} \mu\right) \over \sqrt \pi \ |\Gamma(\mu)| \ 2^{(1- 2^{1-n})/2}} $$

$$\times \  \bigg[ \sinh \left(2^{n-1} \pi\tau \right)\bigg]^{- 2^{-n}} x^{(1- 2^{-n} )/2- {\rm Re}\mu}, \quad n \in \mathbb{N}, \ x, \tau > 0.\eqno(1.11) $$
Finally in this section,  the case of the Olevskii or Fourier-Jacobi transform [1] with the Gauss hypergeometric function [3], Vol. III in the kernel can be studied, appealing to the following representation (see [3], Vol. II, Entry 2.16.21.1)

$${}_2F_1\left( {\mu+\nu\over 2} + i\tau, \   {\mu+\nu\over 2} - i\tau;  \nu +1;  - x^2\right) $$

$$= {2^{2-\mu} x^{-\nu}\  \Gamma(\nu+1) \over |\Gamma ((\mu+\nu)/2 +i\tau)|^2} \int_0^\infty y^{\mu-1} J_\nu(xy) K_{2i\tau}(y) dy,\quad x, \tau, \ \mu+\nu >0.\eqno(1.12)$$ 
Therefore,  the Cauchy-Shcwarz inequality and Entries 2.12.31.1, 2.16.33.2,  in [3], Vol. II yield

$$ \left| {}_2F_1\left( {\mu+\nu\over 2} + i\tau, \   {\mu+\nu\over 2} - i\tau;  \nu +1;  - x^2\right) \right| $$

$$\le  {2^{2-\mu} x^{-\nu-\mu/2}\  \Gamma(\nu+1) \over |\Gamma ((\mu+\nu)/2 +i\tau)|^2} \left( \int_0^\infty y^{\mu-1} J^2_\nu(y) dy\right)^{1/2} \left( \int_0^\infty y^{\mu-1} K^2_{2i\tau}(y) dy\right)^{1/2}$$ 

$$=   2^{1/2-\mu}   \Gamma(\nu+1) \bigg[ {\Gamma (\mu/2)\Gamma ((1-\mu)/2) \Gamma(\mu/2+\nu)\over   \Gamma(1-\mu/2) \Gamma(1-\mu/2+\nu) \Gamma((1+\mu)/2)}\bigg]^{1/2}$$

$$\times   { |\Gamma (\mu/2 + 2i\tau)| \over   |\Gamma ((\mu+\nu)/2 +i\tau)|^2}\ x^{-\nu-\mu/2},$$
i.e. after simplification

$$ \left| {}_2F_1\left( {\mu+\nu\over 2} + i\tau, \   {\mu+\nu\over 2} - i\tau; \  \nu +1;  - x^2\right) \right| $$

$$\le   2^{1/2-\mu} \  {\Gamma \left({\mu/ 2}\right) \Gamma(\nu+1) \over  \Gamma((1+\mu)/2)}   \bigg[ {\tan(\pi\mu/2)\   \Gamma(\mu/2+\nu)\over  \Gamma(1-\mu/2+\nu)}\bigg]^{1/2}$$

$$\times   { |\Gamma (\mu/2 + 2i\tau)| \over   |\Gamma ((\mu+\nu)/2 +i\tau)|^2}\ x^{-\nu-\mu/2},\quad x,\ \tau >0,\ \nu > - {\mu\over 2},\ 0 <\mu < 1.\eqno(1.13)$$

\section {Uniform asymptotic expansions}

In this section we establish uniform asymptotic expansions for the index kernels above similar to Theorem 1 for the Kontorovich-Lebedev kernel $K_{i\tau}(x)$, containing explicit remainders with the corresponding estimates.   As  is known, the modified Bessel function  $K_{i\tau}(x)$ has the representation 
$$K_{i\tau}(x)= {\pi\over 2i\sinh(\pi\tau)} \left[ I_{-i\tau}(x)- I_{i\tau}(x)\right],\quad x,  \tau  >0,\eqno(2.1)$$ 
where the modified Bessel function of the first kind $I_{\nu}(x)$ is defined by the series 

$$I_{\nu}(x)= \sum_{k=0}^\infty {(x/2)^{2k+\nu} \over k! \Gamma(k+\nu+1)}.\eqno(2.2)$$
Meanwhile, the kernel of the Lebedev transform, which is the  square of the modified Bessel function $K^2_{i\tau}(x)$,  is written  in terms of the following hypergeometric ${}_1F_2$-functions (see [1, formula (1.166)])

$$K^2_{i\tau}(x) =  {\pi \over 2 \tau \sinh(\pi\tau)} {}_1F_2 \left({1\over 2}; 1+i\tau,\ 1-i\tau; x^2 \right) $$

$$+ {1\over 2} {\rm Re} \left[ \left({x\over 2}\right)^{-2i\tau}  \Gamma^2(i\tau)\ {}_1F_2 \left({1\over 2} - i\tau; 1-i\tau,\ 1- 2i\tau; x^2 \right)\right].\eqno(2.3)$$
We treat each hypergeometric function with the use of (2.2) and the beta integral,   interchanging the order of integration and summation by virtue of the absolute and uniform convergence.  Precisely, we deduce

$${}_1F_2 \left({1\over 2}; 1+i\tau,\ 1-i\tau; x^2 \right) =1+  {\sqrt \pi \tau \over \sinh(\pi\tau)} \sum_{k=1}^\infty { \Gamma (1/2+k) \over \Gamma(1+ i\tau+k) \Gamma(1- i\tau+k)} {x^{2k}\over k!}$$

$$=  1+  {2 (2x)^{2i\tau} \Gamma(i\tau) \over \sqrt \pi  \Gamma(1/2+ i\tau)}  \int_0^1\int_0^{2xu} (1-u^2)^{i\tau-1/2} u^{2i\tau} (4x^2u^2 -v^2)^{-i\tau} I_1(v)  dv du,\eqno(2.4)$$

$${}_1F_2 \left({1\over 2} - i\tau; 1-i\tau,\ 1- 2i\tau; x^2 \right) = 1+  { 2^{-2i\tau} \Gamma(1- i\tau) \over \sqrt \pi \ \Gamma (1/2-i\tau)} $$

$$\times \sum_{k=1}^\infty {B\left(1/2 - i\tau,  1/2 - i\tau +k\right) B\left(1 - i\tau,  k\right) \ x^{2k}\over k! (k-1)!}= 1+  { 2 x^{2i\tau}  \Gamma(1- i\tau) \over \sqrt \pi \ \Gamma (1/2-i\tau)} $$

$$\times  \int_0^1\int_0^{2xu} (1-u^2)^{-i\tau-1/2} (4x^2u^2 -v^2)^{-i\tau} I_1(v)   dv du,\eqno(2.5)$$
where $B(\alpha, \beta)$ is the Euler beta function. Consequently, combining with (2.3), (2.4), (2.5) and integrating by parts in the corresponding inner integrals, we establish the  representation for the square of the modified Bessel function

$$ K^2_{i\tau}(x) =  { \pi  \over 2 \tau \sinh(\pi\tau)} \bigg[ 1+  {\rm Re} \left[ \left({x\over 2}\right)^{-2i\tau}  {\Gamma (i\tau)\over \Gamma(-i\tau)} \left[ 1+ {1\over 2(1-i\tau)}\right] \right] - {x^2i\over \tau(1+\tau^2)}\bigg.$$

$$+  {2 (2x)^{2i\tau} \Gamma(i\tau) \over \sqrt \pi  \Gamma(1/2+ i\tau) (1-i\tau) }  \int_0^1\int_0^{2xu} (1-u^2)^{i\tau-1/2} u^{2i\tau}  (4x^2u^2 -v^2)^{1-i\tau} I_2(v)  {dv du\over v} $$

$$  -  {1\over \pi} \cosh(\pi\tau) {\rm Re} \left[ {  \Gamma(1+2i\tau) \over 1-i\tau} \int_0^1\int_0^{2xu} (1-u^2)^{-i\tau-1/2} \right.$$

$$\left.\bigg. \times (4x^2u^2 -v^2)^{1-i\tau} I_2(v)   {dv du\over v} \right] \bigg].\eqno(2.6)$$
This is a key ingredient which gives the asymptotic expansion by the index at infinity for the  square of the modified Bessel function. Indeed, in the same manner as in [2], Theorem 2.1 one proves

{\bf Theorem 2}. {\it Let  $ x \in \left(0, X\right],\ X >0$. Then the square of the modified Bessel function $K^2_{i\tau}(x)$  has the following asymptotic expansion

$$ K^2_{i\tau}(x) =  { \pi  \over 2 \tau \sinh(\pi\tau)} \bigg[ 1+  \sin \left( 2\tau \log \left({2\tau\over e x} \right)  \right) + R(\tau) \bigg], \quad \tau \to +\infty,\eqno(2.7)$$
where the remainder term $R(\tau)$ is given explicitly

$$R(\tau) = - \left(1+ e^{-2\pi \tau} \right)  {\rm Re} \left[ i \exp\left( 2 i\tau \log \left({2\tau\over e x} \right)  \right) r(i\tau) \bigg[ 1+ {r(\tau)\over 2} \bigg] {3-2i\tau\over 1-i\tau} \right] \bigg.$$

$$+  {2 (2x)^{2i\tau} \Gamma(i\tau) \over \sqrt \pi  \Gamma(1/2+ i\tau) (1-i\tau) }  \int_0^1\int_0^{2xu} (1-u^2)^{i\tau-1/2} u^{2i\tau}  (4x^2u^2 -v^2)^{1-i\tau} I_2(v)  {dv du\over v} $$

$$- {x^2i\over \tau(1+\tau^2)}  -   {1\over \pi} \cosh(\pi\tau) {\rm Re} \left[ {  \Gamma(1+2i\tau) \over 1-i\tau} \int_0^1\int_0^{2xu} (1-u^2)^{-i\tau-1/2}\right.$$

$$\left. \bigg. \times  (4x^2u^2 -v^2)^{1-i\tau} I_2(v)   {dv du\over v} \right] \eqno(2.8)$$
and $r(i\tau)$ is defined by $(1.4)$. Besides,  $R(\tau),\  \tau \ge \tau_0 >0$  has  the upper bound in the form }

$$\left|R(\tau)\right| \le  {4X  (\coth(\pi\tau_0) )^{1/2} \over \sqrt {\pi \tau} } I_1(2X)  $$

$$ +  { 1 \over \tau}  \bigg[ {8X  (\coth(\pi\tau_0) )^{1/2} \over \sqrt {\pi \tau_0} } I_1(2X) + e^{1/(6\tau_0)} + 2X^2 +   { e^{1/(3\tau_0)} \over 12 \tau_0} \bigg].\eqno(2.9)$$

Further, one treats the kernel (1.8), employing the representation (cf. [1, formula (1.169)] )

$$[I_{i\tau}(x)+  I_{-i\tau}(x)] K_{i\tau} (x) =  {\rm Re} \left[ \left({x\over 2}\right)^{2i\tau} { \Gamma(-i\tau)\over \Gamma(1+i\tau)}\right.$$

$$\left. \times   {}_1F_2 \left({1\over 2} + i\tau; 1+ i\tau,\ 1+ 2i\tau; x^2 \right)\right],\quad x, \tau > 0.\eqno(2.10)$$
Hence,  recalling (2.5), (2.6),  we find 

$$[I_{i\tau}(x)+  I_{-i\tau}(x)] K_{i\tau} (x) =  {\rm Re} \left[  \left({x\over 2}\right)^{2i\tau} { \Gamma(-i\tau)\over \Gamma(1+i\tau)} \left[ 1+  {  x^{2} \over  1+i\tau}\right] \right.$$ 

$$\left. + {  2^{-2i\tau}  \Gamma(- i\tau) \over   \sqrt \pi \ \Gamma (1/2+i\tau)(1+i\tau)}    \int_0^1\int_0^{2xu} (1-u^2)^{i\tau-1/2} (4x^2u^2 -v^2)^{1+i\tau} I_2(v)   {dv\over v}  du\right] $$

$$= - {\sinh(\pi\tau)\over \pi}  {\rm Im} \left[  \left({x\over 2}\right)^{2i\tau}  \Gamma^2(i\tau)\right] \left[ 1+ {  x^{2} \over  1+\tau^2}\right] +  {  x^{2} \tau \ \sinh(\pi\tau) \over  \pi(1+\tau^2)} \  {\rm Im} \left[  \left({x\over 2}\right)^{-2i\tau}  \Gamma^2(i\tau)\right]$$

$$+  {\rm Re} \left[  {  2^{-2i\tau}  \Gamma(- i\tau) \over   \sqrt \pi \ \Gamma (1/2+i\tau) (1+i\tau)}  \int_0^1\int_0^{2xu} (1-u^2)^{i\tau-1/2}\right.$$

$$\left. \times  (4x^2u^2 -v^2)^{1+i\tau} I_2(v)  { dv \over v} du\right].\eqno(2.11)$$
Then, appealing to   [2, formula  (2.12)], equality (2.11)  yields 

$$[I_{i\tau}(x)+  I_{-i\tau}(x)] K_{i\tau} (x) =  {1\over \tau} \left(1- e^{-2\pi\tau}\right) \cos\left( 2\tau \log\left({x\tau\over 2 e}\right) \right) \left(1+r(\tau)\right)^2$$

$$\times   \left[ 1+  {  x^{2} (1-\tau) \over  1+ \tau^2}\right] +  {\rm Re} \left[  {  2^{-2i\tau}  \Gamma(- i\tau) \over   \sqrt \pi \ \Gamma (1/2+i\tau)(1+i\tau)} \right.$$

$$\left. \times \int_0^1\int_0^{2xu} (1-u^2)^{i\tau-1/2} (4x^2u^2 -v^2)^{1+i\tau} {I_2(v)\over v}    dv du\right].\eqno(2.12)$$

{\bf Theorem 3}. {\it Let  $ x \in \left(0, X\right],\ X >0$. Then the kernel $(1.8)$   has the following asymptotic expansion

$$ [I_{i\tau}(x)+  I_{-i\tau}(x)] K_{i\tau} (x) =  { 1  \over  \tau} \bigg[ \cos\left( 2\tau \log\left({x\tau\over 2 e}\right) \right)  + R(\tau) \bigg], \quad \tau \to +\infty,\eqno(2.13)$$
where the remainder term $R(\tau)$ is given by the formula 

$$R(\tau) =  \cos\left( 2\tau \log\left({x\tau\over 2 e}\right) \right) \left[ (1-  e^{-2\pi\tau}  )\ \left[  r(i\tau)\ (2+r(i\tau) )\left( 1+ {  x^{2} (1-\tau) \over  1+ \tau^2} \right)\right.\right.$$

$$\left.\left. +  {  x^{2} (1-\tau) \over  1+ \tau^2} \right] -  e^{-2\pi\tau}  \right] + {\coth (\pi\tau)\over \sqrt \pi}  {\rm Re} \left[  {  2^{-2i\tau}  i \tau \ \Gamma(1/2 - i\tau) \over   \Gamma (2+i\tau)}\right.$$

$$\left. \times \int_0^1\int_0^{2xu} (1-u^2)^{i\tau-1/2} (4x^2u^2 -v^2)^{1+i\tau} {I_2(v)\over v}    dv du\right]. \eqno(2.14)$$
Furthermore,  $R(\tau),\  \tau \ge \tau_0 >0$  has  the upper bound}

$$\left|R(\tau)\right| <  {1\over \sqrt \tau} \bigg[  {X \sqrt\pi  \over  \sqrt{  \tanh (\pi\tau_0)}} \  I_1(2 X) +  {1-  e^{-2\pi\tau_0} \over 6 \sqrt \tau } \   e^{1/(3\tau_0)} \left( 1+  {X^{2}\over \tau_0} \left(1+  {1 \over \tau_0} \right)\right)$$

$$ +  {  X^{2}\over \sqrt \tau} \  \left(1+  {1 \over \tau_0} \right) + \sqrt \tau \  e^{-2\pi\tau}\bigg].\eqno(2.15)$$

\begin{proof} In fact,  writing  (2.12) in the form

$$[I_{i\tau}(x)+  I_{-i\tau}(x)] K_{i\tau} (x) =  {1\over \tau} \left[  \cos\left( 2\tau \log\left({x\tau\over 2 e}\right) \right)\left( 1+ r(i\tau)(2+r(i\tau) )\left( 1+ {  x^{2} (1-\tau) \over  1+ \tau^2} \right)\right.\right.$$

$$\left. \left. +  {  x^{2} (1-\tau) \over  1+ \tau^2} \right) + {\coth (\pi\tau)\over \sqrt \pi}  {\rm Re} \left[  {  2^{-2i\tau}  i \tau \ \Gamma(1/2 - i\tau) \over   \Gamma (2+i\tau)} \right.\right.$$

$$\left. \times \int_0^1\int_0^{2xu} (1-u^2)^{i\tau-1/2} (4x^2u^2 -v^2)^{1+i\tau} {I_2(v)\over v}    dv du\right]$$

$$\left. -  e^{-2\pi\tau}  \cos\left( 2\tau \log\left({x\tau\over 2 e}\right) \right) \left(1+r(i\tau)\right)^2 \left[ 1+  {  x^{2} (1-\tau) \over  1+ \tau^2}\right] \right],$$
we easily end up with (2.13), (2.14).  Finally, the desired estimate (2.15) one gets, using the bound for $r(t)$ (cf. [2])

$$\left|r(i\tau)\right| \le e^{1/(6\tau)} -1<  {e^{1/(6\tau)} \over 6\tau}.\eqno(2.16)$$ 
Indeed, we have for $\tau \ge \tau_0 >0$

$$ | R(\tau) | <  {1-  e^{-2\pi\tau_0} \over 6\tau } \   e^{1/(3\tau_0)} \left( 1+  {X^{2}\over \tau_0} \left(1+  {1 \over \tau_0} \right)\right)$$

$$ +  {  X^{2}\over \tau} \  \left(1+  {1 \over \tau_0} \right) +  e^{-2\pi\tau} + 4 X^2\  {\coth (\pi\tau)\over \sqrt \pi}   \left|{ \Gamma(1/2 - i\tau) \over   \Gamma (1+i\tau)}\right| $$

$$ \times \int_0^1 (1-u^2)^{-1/2}  du \int_0^{2X}  {I_2(v)\over v}   dv$$

$$= {1-  e^{-2\pi\tau_0} \over 6\tau } \   e^{1/(3\tau_0)} \left( 1+  {X^{2}\over \tau_0} \left(1+  {1 \over \tau_0} \right)\right)$$

$$ +  {  X^{2}\over \tau} \  \left(1+  {1 \over \tau_0} \right) +  e^{-2\pi\tau} + {X \sqrt\pi  \over  \sqrt{ \tau \tanh (\pi\tau_0)}} \  I_1(2 X)$$
which yields  (2.15) and completes the proof. 

\end{proof}

Let us treat the kernel for the index Whittaker transform which is given by formula (1.171) in [1] in terms of the confluent hypergeometric function ${}_1F_1$ and can be rewritten due to Entry 7.11.1.2 in [3], Vol. III in the form

$$   e^{-x/2}\  W_{\rho, i\tau} (x) = 2 {\rm Re} \left[  { \Gamma(-2 i\tau)\ x^{i\tau+1/2} \over \Gamma(1/2-\rho-i\tau)}  \   {}_1F_1 \left({1\over 2} +\rho + i\tau; \ 1+ 2i\tau;  \  - x \right)\right],\eqno(2.17)$$
where $\  x, \tau > 0,\ \rho \in \mathbb{R}$.  In the meantime, employing the known integral representation for the confluent hypergeometric function, we have with simple substitutions for $|\rho| <  1/2$

$$ {}_1F_1 \left({1\over 2} +\rho + i\tau; \ 1+ 2i\tau;  \  - x \right) = { \Gamma(1+2 i\tau)\over  \Gamma(1/2+\rho+i\tau) \Gamma(1/2-\rho+i\tau)} $$

$$\times \int_0^1 e^{-xy} y^{i\tau+\rho-1/2}(1-y)^{i\tau-\rho-1/2} dy =  { \Gamma(1+2 i\tau)\ 4^{-i\tau} \ e^{-x/2} \over  \Gamma(1/2+\rho+i\tau) \Gamma(1/2-\rho+i\tau)}  $$

$$\times \bigg[ \int_0^1 e^{xy/2} (1-y)^{i\tau+\rho-1/2}(1+y)^{i\tau-\rho-1/2} dy +  \int_0^1 e^{-xy/2} (1-y)^{i\tau-\rho-1/2}(1+y)^{i\tau+\rho-1/2} dy\bigg] $$

$$ =  { \Gamma(1+2 i\tau)\ 4^{-i\tau} \ e^{-x/2} \over  \Gamma(1/2+\rho+i\tau) \Gamma(1/2-\rho+i\tau)} \int_{-1}^1 e^{-xy/2} (1-y)^{i\tau-\rho-1/2}(1+y)^{i\tau+\rho-1/2} dy $$

$$ =  e^{-x/2} \bigg[ 1+  { \Gamma(1+2 i\tau)\ 4^{-i\tau} \over  \Gamma(1/2+\rho+i\tau) \Gamma(1/2-\rho+i\tau)} \int_{-1}^1 \left(e^{-xy/2} -1\right) (1-y)^{i\tau-\rho-1/2}(1+y)^{i\tau+\rho-1/2} dy \bigg].$$
Then, invoking Entry 2.2.6.1 in [3], Vol. I, we find the equality

$$ {}_1F_1 \left({1\over 2} +\rho + i\tau; \ 1+ 2i\tau;  \  - x \right)  = e^{-x/2} \bigg[ 1+  \sum_{k=1}^\infty {(x/2)^k\over k!} \bigg.$$

$$\bigg. \times \ {}_2F_1\left(-k;\  i\tau+\rho+ {1\over 2}; \ 1+2i\tau;\ 2 \right) \bigg].\eqno(2.18)$$
Meanwhile,

$${}_2F_1\left(-k;\  i\tau+\rho+ {1\over 2}; \ 1+2i\tau;\ 2 \right) = \sum_{m=0}^k {2^m\over m!} \ {(-k)_m\  (i\tau+\rho + 1/2)_m\over (1+2i\tau)_m}$$

$$= \sum_{m=0}^k (-1)^m \binom{k}{m}  \prod_{j=1}^{m} \left( 1+ {2\rho+j-1\over 2i\tau+j} \right).$$
Writing the latter product in the form

$$\prod_{j=1}^{m} \left( 1+ {2\rho+j-1\over 2i\tau+j} \right) = 1+ \sum_{r=1}^m\  \sum_{1\le j_1<j_2<\dots < j_r \le m} \ \prod_{\mu=1}^r  {2\rho+ j_\mu -1\over 2i\tau+j_\mu},$$
it implies 

$${}_2F_1\left(-k;\  i\tau+\rho+ {1\over 2}; \ 1+2i\tau;\ 2 \right) = \sum_{m=1}^k (-1)^m \binom{k}{m} \sum_{r=1}^m\  \sum_{1\le j_1<j_2<\dots < j_r \le m} \ \prod_{\mu=1}^r  {2\rho+ j_\mu -1\over 2i\tau+j_\mu}$$

$$=  \sum_{r=1}^k\  \sum_{m=r}^k (-1)^m \binom{k}{m}  \sum_{1\le j_1<j_2<\dots < j_r \le m} \ \prod_{\mu=1}^r  {2\rho+ j_\mu -1\over 2i\tau+j_\mu},$$
and 

$$\left| {}_2F_1\left(-k;\  i\tau+\rho+ {1\over 2}; \ 1+2i\tau;\ 2 \right) \right| \le \sum_{r=1}^k\  \sum_{m=r}^k  \binom{k}{m}  \sum_{1\le j_1<j_2<\dots < j_r \le m} \ \prod_{\mu=1}^r  \left|{2\rho+ j_\mu -1\over 2i\tau+j_\mu}\right| $$

$$ \le 2^k \sum_{r=1}^k\  {1\over (2\tau)^r }  \sum_{1\le j_1<j_2<\dots < j_r \le k} \ \prod_{\mu=1}^r  (2|\rho|+ j_\mu +1) \le {2^{k-1}\over \tau} \left[ \prod_{r=1}^{k} \left( 2|\rho| + r +1\right)- 1 \right]$$

$$ =  {2^{k-1}\over \tau} \left[\ (2(|\rho|+1))_k -1\ \right].$$
Consequently, one has the following estimate from (2.18)

$$\left|  {}_1F_1 \left({1\over 2} +\rho + i\tau; \ 1+ 2i\tau;  \  - x \right)\right| \le e^{-x/2} \bigg[ 1+  {1\over 2\tau } \bigg[ \sum_{k=1}^\infty {x^k  (2(|\rho|+1))_k \over k!}  +1 -  e^x \bigg] \bigg]$$

$$=  e^{-x/2} \bigg[ 1+  {1\over 2\tau } \bigg[ (1-x)^{- 2(1+|\rho|)}  -  e^x \bigg] \bigg], \quad \tau > 0,\ 0 < x < 1.\eqno(2.19)$$
Meanwhile, having in mind   the representation for the gamma function (cf. [4])

$$\Gamma(z)  =   \sqrt{2\pi}\  \exp \left( \left(z- {1\over 2} \right) \log z - z \right) \left[1 + r(z)\right],\quad  |\arg z | \le {\pi\over 2},\eqno(2.20)$$
where 

$$r(z) = \exp\left( \int_0^\infty e^{- z t} \left[ {1\over 2}- {1\over t} + {1\over e^t-1}\right] {dt\over t} \right) -1\eqno(2.21)$$
and (see (2.16))

$$\left|r(z)\right| \le e^{1/(6 |z|)} -1\le {e^{1/(6 |z|)}\over 6|z|} ,\eqno(2.22)$$
we see that it generalizes (1.4).   Moreover,  it yields

$$ \Gamma\left({1\over 2} +\rho+ i\tau\right) =  \sqrt{2\pi}\  \exp \left( \left(\rho +i\tau \right) \log \left( {1\over 2} +\rho +i\tau \right)  - {1\over 2} -\rho-i\tau \right) $$

$$\times \left[1 + r\left({1\over 2} +\rho+i\tau\right)\right],\eqno(2.23)$$

$$ \Gamma\left(- 2 i\tau\right) =  \sqrt{2\pi}\  \exp \left( \left(-2 i\tau- {1\over 2}  \right) \log \left(-2 i\tau \right)  + 2i\tau \right) \left[1 + r\left( -2 i\tau\right)\right].\eqno(2.24)$$
Then, returning to (2.17) and taking into account (2.18), (2.23), (2.24),  we write it in the form 

$$    W_{\rho, i\tau} (x) =  \sqrt {2 x}  \ e^{- \pi\tau/2 } \ \tau^{\rho-1/2}  \exp \left( \tau \arctan \left({1+2\rho \over 2\tau  } \right) - {1\over 2} -\rho +  {\rho \over 2} \  \log \left(1+  {(1 + 2\rho)^2 \over 4\tau^2}  \right)\right) $$

$$\times  \ {\rm Re} \left[ \left( e^{\pi i \rho - 2\pi \tau} +  e^{- \pi i \rho} \right)  \exp \left(  i\tau  \  \log \bigg[  {e x \over 4\tau} \ \left(1+  {(1 + 2\rho) \over 4\tau^2 }  \right)^{1/2} \bigg]  - i \rho \arctan \left({1+2\rho \over 2\tau  } \right) \right.\right.$$

$$\left.  + {\pi i\over 2} \left(\rho + {1\over 2} \right) \right)    \left[1 + r\left({1\over 2} +\rho+i\tau\right)\right]  \left[1 + r\left( -2 i\tau\right)\right] $$

$$\left. \times \bigg[ 1+  \sum_{k=1}^\infty {(x/2)^k\over k!}  \ {}_2F_1\left(-k;\  i\tau+\rho+ {1\over 2}; \ 1+2i\tau;\ 2 \right) \bigg] \right].\eqno(2.25)$$

 We summarize the above results as the following theorem.

{\bf Theorem 4}. {\it Let  $ x \in \left(0, x_0\right],\  0 < x_0 < 1,\  \rho \in \mathbb{R},\ |\rho| <  1/2$. Then  kernel $(2.17)$   has the  asymptotic expansion

$$   W_{\rho, i\tau} (x) =   \sqrt {2 x}  \ e^{- \pi\tau/2 } \ \tau^{\rho-1/2}  \exp \left( \tau \arctan \left({1+2\rho \over 2\tau  } \right) - {1\over 2} -\rho +  {\rho \over 2} \  \log \left(1+  {(1 + 2\rho)^2 \over 4\tau^2}  \right)\right) $$

$$\times \left[  \cos \left(  \tau  \  \log \bigg[  {e x \over 4\tau} \ \left(1+  {(1 + 2\rho) \over 4\tau^2 }  \right)^{1/2} \bigg]  -  \rho \arctan \left({1+2\rho \over 2\tau  } \right)  - {\pi \over 2} \left(\rho - {1\over 2} \right) \right)\right.$$

$$\left.  + R(\tau)\right],\quad \tau \to +\infty,\eqno(2.26)$$
where the remainder term $R(\tau)$ is given by the formula 
$$R(\tau) =  {\rm Re} \left[ \exp \left(  i\tau  \  \log \bigg[  {e x \over 4\tau} \ \left(1+  {(1 + 2\rho) \over 4\tau^2 }  \right)^{1/2} \bigg]  - i \rho \arctan \left({1+2\rho \over 2\tau  } \right)  - {\pi i\over 2} \left(\rho - {1\over 2} \right) \right) \right.$$

$$\times \left. \bigg[   \left[1 + r\left({1\over 2} +\rho+i\tau\right)\right]  \left[1 + r\left( -2 i\tau\right)\right] \bigg[ 1+  \sum_{k=1}^\infty {(x/2)^k\over k!}  \ {}_2F_1\left(-k;\  i\tau+\rho+ {1\over 2}; \ 1+2i\tau;\ 2 \right) \bigg] -1\bigg] \right]$$

$$+ e^{ - 2\pi \tau}  {\rm Re} \left[ \exp \left(  i\tau  \  \log \bigg[  {e x \over 4\tau} \ \left(1+  {(1 + 2\rho) \over 4\tau^2 }  \right)^{1/2} \bigg]  - i \rho \arctan \left({1+2\rho \over 2\tau  } \right)  + {\pi i\over 2} \left(3 \rho + {1\over 2} \right) \right) \right.$$

$$\times   \left[1 + r\left({1\over 2} +\rho+i\tau\right)\right]  \left[1 + r\left( -2 i\tau\right)\right] $$

$$\left. \times \bigg[ 1+  \sum_{k=1}^\infty {(x/2)^k\over k!}  \ {}_2F_1\left(-k;\  i\tau+\rho+ {1\over 2}; \ 1+2i\tau;\ 2 \right) \bigg] \right].$$
Furthermore,  $R(\tau),\  \tau \ge \tau_0 >0$  has  the upper bound}

$$\left|R(\tau)\right| <   {1\over  \tau} \bigg[  e^{1/ \tau_0} \left( 2 + e^{1/ \tau_0} \right) +   \left( 1 + {e^{1/ \tau_0}\over \tau}  \right)^2 \bigg.$$

$$\bigg. \times \left(  (1-x_0)^{- 2(1+|\rho|)}  -  e^{x_0}  +   e^{ - 2\pi \tau} \left(\tau +  (1-x_0)^{- 2(1+|\rho|)}  -  e^{x_0} \right) \right) \bigg].$$

\begin{proof} The proof is straightforward  from (2.25) and estimates (2.19), (2.22).

\end{proof}

Finally, our goal is to investigate the kernel (1.12) of the Olevskii or Fourier-Jacobi transform.  To do this,  we choose some  $T > 1$, splitting  the integral (1.12) as follows

$${}_2F_1\left( {\mu+\nu\over 2} + i\tau, \   {\mu+\nu\over 2} - i\tau;  \nu +1;  - x^2\right) $$

$$= {2^{2-\mu} x^{-\nu}\  \Gamma(\nu+1) \over |\Gamma ((\mu+\nu)/2 +i\tau)|^2} \bigg[ \int_0^T  y^{\mu-1} J_\nu(xy) K_{2i\tau}(y) dy +  \int_T^\infty   y^{\mu-1} J_\nu(xy)\bigg.$$

$$\bigg.\times  K_{2i\tau}(y) dy\bigg].\eqno(2.27)$$
Let us treat the integral over $[T, \infty)$.  Indeed, we employ equality (1.23) in [2] for Bessel function

$$J_\nu(z)= \left({2\over \pi z}\right)^{1/2} \left( \cos \left( z - {\pi\over 2} \nu - {\pi\over 4} \right) \left( \sum_{n=0}^{N-1} (-1)^n {a_{2n}(\nu)\over z^{2n}} + R_{2N} (z,\nu)\right)\right.$$

$$\left. -  \sin \left( z - {\pi\over 2} \nu - {\pi\over 4} \right) \left( \sum_{n=0}^{M-1} (-1)^n {a_{2n+1}(\nu)\over z^{2n+1}} - R_{2M+1} (z,\nu)\right)\right),$$
where 

$$a_n(\nu)= (-1)^n {\cos(\pi\nu)\over 2^n n! \pi} \Gamma\left(n+{1\over 2}+\nu\right) \Gamma\left(n+{1\over 2}-\nu\right),$$
$$R_N(z,\nu)= (-1)^{[N/2]} \left({2\over \pi }\right)^{1/2}  {\cos(\pi\nu)\over z^N \pi} \int_0^\infty {t^{N-1/2} e^{-t} \over 1+ (t/z)^2} K_\nu(t) dt,$$
provided $|{\rm Re} \nu | < N- 1/2,\  |\arg z| < \pi/ 2$.  Then, assuming  conditions  $|\nu| < 1/2,\  5/4 < \mu < 3/2$, we have 

$$\int_T^\infty   y^{\mu-1} J_\nu(xy) K_{2i\tau}(y) dy =  \left({2\over \pi x}\right)^{1/2} \int_T^\infty   y^{\mu-3/2}  K_{2i\tau}(y) \left( \cos \left( xy - {\pi\over 2} \nu - {\pi\over 4} \right)\right. $$

$$\left.  - {1-\nu^2 \over  8 xy}\ \sin \left( xy  - {\pi\over 2} \nu - {\pi\over 4} \right)  \right)dy$$

$$+  \left({2\over \pi x}\right)^{1/2} \int_T^\infty   y^{\mu-3/2}  K_{2i\tau}(y) \cos \left( xy - {\pi\over 2} \nu - {\pi\over 4} \right) \left( \sum_{n=1}^{N-1} (-1)^n {a_{2n}(\nu)\over (xy)^{2n}} + R_{2N} (xy,\nu)\right) dy $$

$$-  \left({2\over \pi x}\right)^{1/2} \int_T^\infty   y^{\mu-3/2}  K_{2i\tau}(y) \sin \left( xy - {\pi\over 2} \nu - {\pi\over 4} \right)$$

$$\times  \left( \sum_{n=1}^{M-1} (-1)^n {a_{2n+1}(\nu)\over (xy)^{2n+1}} - R_{2M+1} (xy,\nu)\right) dy.\eqno(2.28)$$
Hence, in turn,

$$\left({2\over \pi x}\right)^{1/2} \int_T^\infty   y^{\mu-3/2}  K_{2i\tau}(y) \left( \cos \left( xy - {\pi\over 2} \nu - {\pi\over 4} \right)  - {1-\nu^2 \over  8 xy}\ \sin \left( xy  - {\pi\over 2} \nu - {\pi\over 4} \right)  \right)dy$$

$$= \left({2\over \pi x}\right)^{1/2} \bigg[ \cos \left(  {\pi\over 2} \left(\nu + {1\over 2} \right) \right)+  {1-\nu^2 \over  8 x}   \sin \left(  {\pi\over 2} \left(\nu + {1\over 2} \right) \right)\bigg] \int_T^\infty   y^{\mu-3/2}  K_{2i\tau}(y)  \cos ( xy ) dy $$

$$+ \left({2\over \pi x}\right)^{1/2} \bigg[ \sin \left(  {\pi\over 2} \left(\nu + {1\over 2} \right) \right)-  {1-\nu^2 \over  8 x}   \cos \left(  {\pi\over 2} \left(\nu + {1\over 2} \right) \right)\bigg] \int_T^\infty   y^{\mu-5/2}  K_{2i\tau}(y)  \sin ( xy ) dy.$$
The latter two integrals can be treated with the aid of Entry 2.16.14.3 in [3], Vol. II.  In fact, we find

$$\int_T^\infty   y^{\mu-3/2}  K_{2i\tau}(y)  \cos ( xy ) dy = {\pi\over 4}\  {\Gamma(\mu-1/2+2i\tau)\over (x^2+1)^{(\mu-1/2)/2}\  \sin(\pi( (\mu-1/2)/2 - i\tau))} $$

$$\times \bigg[ P^{-2i\tau}_{\mu-3/2} \left({x\over \sqrt{x^2+1}}\right) +  P^{-2i\tau}_{\mu-3/2} \left(- {x\over \sqrt{x^2+1}}\right) \bigg]$$

$$- \int_0^T   y^{\mu-3/2}  K_{2i\tau}(y)  \cos ( xy ) dy,\eqno(2.29)$$

$$\int_T^\infty   y^{\mu-5/2}  K_{2i\tau}(y)  \sin ( xy ) dy = {\pi\over 4}\  {\Gamma(\mu-3/2+2i\tau)\over (x^2+1)^{(\mu-3/2)/2}\  \cos(\pi((\mu-3/2)/2 - i\tau))} $$

$$\times \bigg[ P^{-2i\tau}_{\mu-5/2} \left(- {x\over \sqrt{x^2+1}}\right) -  P^{-2i\tau}_{\mu-5/2} \left( {x\over \sqrt{x^2+1}}\right) \bigg]$$

$$- \int_0^T   y^{\mu-5/2}  K_{2i\tau}(y)  \sin ( xy ) dy.\eqno(2.30)$$
Meanwhile, appealing to the estimate  (1.27) in [2], we observe the following Lebedev type inequalities

$$ \left| K_{i\tau}(x)\right| \le A {(\tau x)^{-1/4}\over \sqrt{\sinh(\pi\tau)}},\quad \tau > 0, \  x \in (0,T],\eqno(2.31)$$
 
 $$ \left| K_{i\tau}(x)\right| \le B {(x/\tau )^{1/4}\over \sqrt{\sinh(\pi\tau)}},\quad \tau > 0, \  x \in [T,\infty).\eqno(2.32)$$
   where $A, B$ are positive absolute constants.  Consequently,  we have the inequalities
  
  $$\left| \int_0^T   y^{\mu-3/2}  K_{2i\tau}(y)  \cos ( xy ) dy \right| \le A \   { T^{\mu-3/4}\  (2\tau)^{-1/4}\over (\mu-3/4) \sqrt{\sinh(2\pi\tau)}},$$

 $$\left| \int_0^T   y^{\mu-5/2}  K_{2i\tau}(y)  \sin ( xy ) dy \right| \le A \   { T^{\mu+1/4}\  (2\tau)^{-1/4}\over (\mu-3/4) \sqrt{\sinh(2\pi\tau)}}.$$
 Moreover, an elementary inequality for Bessel function (cf. [2])  $|J_\nu(x)| \le C_\nu x^{-1/2},  x >0,\ C_\nu >0$ is a constant and (2.31) presume the estimate

$$\left|\int_0^T  y^{\mu-1} J_\nu(xy) K_{2i\tau}(y) dy\right| \le \int_0^T  y^{\mu-1} \left| J_\nu(xy) K_{2i\tau}(y) \right| dy$$

$$\le  A \ C_\nu {(2\tau)^{-1/4} x^{-1/2} T^{\mu-3/4} \over (\mu- 3/4) \sqrt{ \sinh(2\pi\tau)}}.$$ 
Analogously, with the use of (2.32) we find

$$ \left| \int_T^\infty   y^{\mu-3/2}  K_{2i\tau}(y) \cos \left( xy - {\pi\over 2} \nu - {\pi\over 4} \right) \left( \sum_{n=1}^{N-1} (-1)^n {a_{2n}(\nu)\over (xy)^{2n}} + R_{2N} (xy,\nu)\right) dy\right| $$

$$\le  \int_T^\infty   y^{\mu-3/2}  \left| K_{2i\tau}(y)\right| \left( \sum_{n=1}^{N-1} {|a_{2n}(\nu)| \over (xy)^{2n}}  +   {1\over (xy)^{2N}} \int_0^\infty  t^{2N-1/2} e^{-t} K_\nu(t) dt \right) dy$$

$$\le A_{N,\nu,\mu}\  {\tau^{-1/4}\over \sqrt{\sinh(2\pi\tau)}},\quad  x\ge x_0 >0,\ \tau \ge \tau_0 > 0,$$

$$ \left| \int_T^\infty   y^{\mu-3/2}  K_{2i\tau}(y) \sin \left( xy - {\pi\over 2} \nu - {\pi\over 4} \right) \left( \sum_{n=1}^{M} (-1)^n {a_{2n}(\nu)\over (xy)^{2n}} - R_{2M+1} (xy,\nu)\right) dy\right| $$

$$\le  \int_T^\infty   y^{\mu-3/2}  \left| K_{2i\tau}(y)\right| \left( \sum_{n=1}^{M-1} {|a_{2n+1}(\nu)| \over (xy)^{2n+1}}  +   {1\over (xy)^{2M+1}} \int_0^\infty  t^{2M+1/2} e^{-t} K_\nu(t) dt \right) dy$$

$$\le B_{N,\nu,\mu}\  {\tau^{-1/4}\over \sqrt{\sinh(2\pi\tau)}},\quad  x\ge x_0 >0,\ \tau \ge \tau_0 > 0,$$
where $A_{N,\nu,\mu},\  B_{N,\nu,\mu}$ are positive constants and  $ 5/4 < \mu < 3/2,\ |\nu| < 1/2$.

 In the meantime, Entry 2.17.1.6 in [3], Vol. III, differential properties of the  associated Legendre functions and integration by parts multiple times suggest  their  representations as follows  
 
 $$  P^{-2i\tau}_{\bracenom{\mu-3/2} {\mu-5/2}} \left(\pm {x\over \sqrt{x^2+1}}\right) = {  (x^2+1)^{i\tau} \over \Gamma(1+2i\tau)} \left[ \sum_{m=0}^N {(-1)^m\over 2^m\  m! \ (2i\tau+1)_m}  \right.$$
 
 $$\left. \times \Gamma\left( m+1+ \bracenom{\mu-3/2} {\mu-5/2} \right) \left[ \Gamma\left( 1-m+ \bracenom{\mu-3/2} {\mu-5/2} \right) \right]^{-1} \left(1\mp {x\over \sqrt{x^2+1}}\right)^{2i\tau+m} \right.$$
 
 $$\left. + {1\over (2i\tau+1)_N} \int_{\pm x/(x^2+1)^{1/2}}^1 {P^{N+1}_{\bracenom{\mu-3/2} {\mu-5/2}} (t)\over (1-t^2)^{(N+1)/2}}  \left(t\mp {x\over \sqrt{x^2+1}}\right)^{2i\tau} dt \right],\quad N \in \mathbb{N}_0. $$
Hence, accordingly,

$$P^{-2i\tau}_{\mu-3/2} \left({x\over \sqrt{x^2+1}}\right) +  P^{-2i\tau}_{\mu-3/2} \left(- {x\over \sqrt{x^2+1}}\right)$$

$$= { 1\over \Gamma(1+2i\tau)} \left[ 2\sum_{m=0}^N {(-1)^m \cosh \left( (2i\tau+m) \log\left(x+ (x^2+1)^{1/2} \right) \right) \over 2^m\  m! \  (x^2+1)^{m/2} \ (2i\tau+1)_m}  {\Gamma\left( m+ \mu-1/2 \right)\over  \Gamma\left( \mu-m - 1/2 \right)} \right.$$
 
 $$ + {1\over (2i\tau+1)_N} \int_{ x/(x^2+1)^{1/2}}^1 {P^{N+1}_{\mu-3/2} (t)\over (1-t^2)^{(N+1)/2}}  \left(t \sqrt{x^2+1} - x\right)^{2i\tau} dt $$

$$\left. + {1\over (2i\tau+1)_N} \int_{ - x/(x^2+1)^{1/2}}^1 {P^{N+1}_{\mu-3/2} (t)\over (1-t^2)^{(N+1)/2}}  \left(t \sqrt{x^2+1} + x\right)^{2i\tau} dt \right],\eqno(2.33)$$

$$P^{-2i\tau}_{\mu-5/2} \left(- {x\over \sqrt{x^2+1}}\right) -  P^{-2i\tau}_{\mu-5/2} \left( {x\over \sqrt{x^2+1}}\right)$$

$$= { 1\over \Gamma(1+2i\tau)} \left[ 2\sum_{m=0}^N {(-1)^m \sinh \left( (2i\tau+m) \log\left(x+ (x^2+1)^{1/2} \right) \right) \over 2^m\  m! \  (x^2+1)^{m/2} \ (2i\tau+1)_m}  {\Gamma\left( m+ \mu-3/2 \right)\over  \Gamma\left( \mu-m - 3/2 \right)} \right.$$
 
 $$ + {1\over (2i\tau+1)_N} \int_{ - x/(x^2+1)^{1/2}}^1 {P^{N+1}_{\mu-5/2} (t)\over (1-t^2)^{(N+1)/2}}  \left(t \sqrt{x^2+1} + x\right)^{2i\tau} dt $$

$$\left. - {1\over (2i\tau+1)_N} \int_{  x/(x^2+1)^{1/2}}^1 {P^{N+1}_{\mu-5/2} (t)\over (1-t^2)^{(N+1)/2}}  \left(t \sqrt{x^2+1} - x\right)^{2i\tau} dt \right].\eqno(2.34)$$
Consequently, returning to (2.27) and taking into account (2.29), (2.30), (2.33), (2.34) and asymptotic behavior of the associated Legendre functions near the unity, we establish the following asymptotic expansion of the Olevskii kernel (1.12) when $\tau \to + \infty$,  $x \in [x_0, X_0],\ 0 < x_0 < X_0,\  5/4 < \mu < 3/2,  |\nu| < 1/2$

$${}_2F_1\left( {\mu+\nu\over 2} + i\tau, \   {\mu+\nu\over 2} - i\tau;  \nu +1;  - x^2\right) $$

$$= {2^{2-\mu-1/2} \sqrt \pi\ x^{-\nu-1/2}\  (x^2+1)^{(1/2-\mu)/2}\  \Gamma(\nu+1) \ \Gamma(\mu-1/2+2i\tau)   \over |\Gamma ((\mu+\nu)/2 +i\tau)|^2\  \Gamma(1+2i\tau) \sin(\pi((\mu-1/2)/2 - i\tau))}$$

$$\times  \bigg[ \cos \left(  {\pi\over 2}  \left(\nu + {1\over 2} \right) \right)+  {1-\nu^2 \over  8 x}   \sin \left(  {\pi\over 2} \left(\nu + {1\over 2} \right) \right)\bigg] $$

$$\times \cos \left(2\tau \log\left(x+ (x^2+1)^{1/2} \right) \right)+ R_{N,M}(\tau),\quad \tau \to +\infty,\eqno(2.35)$$
where

$$R_{N,M}(\tau) = {2^{2-\mu} x^{-\nu}\  \Gamma(\nu+1) \over |\Gamma ((\mu+\nu)/2 +i\tau)|^2} \bigg[ \int_0^T  y^{\mu-1}  K_{2i\tau}(y)  J_\nu(xy) dy \bigg.$$

$$+\left({2\over \pi x}\right)^{1/2} \bigg[ \cos \left(  {\pi\over 2} \left(\nu + {1\over 2} \right) \right)+  {1-\nu^2 \over  8 x}   \sin \left(  {\pi\over 2} \left(\nu + {1\over 2} \right) \right)\bigg]$$

$$\times \bigg[ {\pi\over 4}\  { (x^2+1)^{(1/2-\mu)/2} \ \Gamma(\mu-1/2+2i\tau)\over \Gamma(1+2i\tau)  \sin(\pi( (\mu-1/2)/2 - i\tau))} \bigg.$$

$$\times  \left[ 2\sum_{m=1}^N {(-1)^m \cosh \left( (2i\tau+m) \log\left(x+ (x^2+1)^{1/2} \right) \right) \over 2^m\  m! \  (x^2+1)^{m/2} \ (2i\tau+1)_m}  {\Gamma\left( m+ \mu-1/2 \right)\over  \Gamma\left( \mu-m - 1/2 \right)} \right.$$
 
 $$ + {1\over (2i\tau+1)_N} \int_{ x/(x^2+1)^{1/2}}^1 {P^{N+1}_{\mu-3/2} (t)\over (1-t^2)^{(N+1)/2}}  \left(t \sqrt{x^2+1} - x\right)^{2i\tau} dt $$

$$\left. + {1\over (2i\tau+1)_N} \int_{ - x/(x^2+1)^{1/2}}^1 {P^{N+1}_{\mu-3/2} (t)\over (1-t^2)^{(N+1)/2}}  \left(t \sqrt{x^2+1} + x\right)^{2i\tau} dt \right]$$

$$\bigg. - \int_0^T  y^{\mu-3/2}  K_{2i\tau}(y) \cos ( xy ) dy\bigg] $$

$$+ \left({2\over \pi x}\right)^{1/2} \bigg[ \sin \left(  {\pi\over 2} \left(\nu + {1\over 2} \right) \right)-  {1-\nu^2 \over  8 x}   \cos \left(  {\pi\over 2} \left(\nu + {1\over 2} \right) \right)\bigg] $$

$$\times \bigg[ {\pi\over 4}\  {(x^2+1)^{(3/2-\mu)/2} \Gamma(\mu-3/2+2i\tau)\over  \Gamma(1+2i\tau) \cos(\pi((\mu-3/2)/2 - i\tau))} \bigg.$$

$$\times \left[ 2\sum_{m=0}^N {(-1)^m \sinh \left( (2i\tau+m) \log\left(x+ (x^2+1)^{1/2} \right) \right) \over 2^m\  m! \  (x^2+1)^{m/2} \ (2i\tau+1)_m}  {\Gamma\left( m+ \mu-3/2 \right)\over  \Gamma\left( \mu-m - 3/2 \right)} \right.$$
 
 $$ + {1\over (2i\tau+1)_N} \int_{ - x/(x^2+1)^{1/2}}^1 {P^{N+1}_{\mu-5/2} (t)\over (1-t^2)^{(N+1)/2}}  \left(t \sqrt{x^2+1} + x\right)^{2i\tau} dt $$

$$\left. - {1\over (2i\tau+1)_N} \int_{  x/(x^2+1)^{1/2}}^1 {P^{N+1}_{\mu-5/2} (t)\over (1-t^2)^{(N+1)/2}}  \left(t \sqrt{x^2+1} - x\right)^{2i\tau} dt \right]$$

$$\bigg. -\int_0^T  y^{\mu-5/2}  K_{2i\tau}(y) \sin ( xy ) dy\bigg]$$

$$+  \left({2\over \pi x}\right)^{1/2} \int_T^\infty   y^{\mu-3/2}  K_{2i\tau}(y) \cos \left( xy - {\pi\over 2} \nu - {\pi\over 4} \right) \left( \sum_{n=1}^{N-1} (-1)^n {a_{2n}(\nu)\over (xy)^{2n}} + R_{2N} (xy,\nu)\right) dy $$

$$-  \left({2\over \pi x}\right)^{1/2} \int_T^\infty   y^{\mu-3/2}  K_{2i\tau}(y) \sin \left( xy - {\pi\over 2} \nu - {\pi\over 4} \right)$$

$$\bigg. \times  \left( \sum_{n=1}^{M-1} (-1)^n {a_{2n+1}(\nu)\over (xy)^{2n+1}} - R_{2M+1} (xy,\nu)\right) dy \bigg].$$
Finally, we note that by virtue of  the Stirling asymptotic formula for the gamma function [4]  the main  term in (2.35) is $O(\tau^{-1/2-\nu}) $ and the remainder $R_{N,M}(\tau) = O(\tau^{3/4-\mu-\nu})$.

\bigskip
\noindent{{\bf Acknowledgments}}
\bigskip

\noindent The work was partially supported by CMUP, which is financed by national funds through FCT (Portugal)  under the project with reference UIDB/00144/2020.

\bigskip
\noindent{{\bf Disclosure statement}}
\bigskip

\noindent No potential conflict  of interest was reported by the author.

\bigskip
\centerline{{\bf References}}
\bigskip
\baselineskip=12pt
\medskip
\begin{enumerate}

\item[{\bf 1.}\ ]  Yakubovich S.  Index Transforms. Singapore:  World Scientific Publishing Company;  1996.

\item[{\bf 2.}\ ]   Yakubovich S.  Upper  bounds and asymptotic expansion for Macdonald's function and the summability of the Kontorovich-Lebedev integrals.  Integral Transforms and Special Functions.  2023; 34, N 10:  721-736.

\item[{\bf 3.}\ ]  Prudnikov AP,  Brychkov  YuA,  Marichev OI. Integrals and series:  Vol. I: Elementary functions. New York:  Gordon and Breach;   1986;   Vol. II:  Special functions. New York: Gordon and Breach;  1986;   Vol. III:  More special functions. New York:   Gordon and Breach; 1990.

\item[{\bf 4.}\ ]  Lebedev NN.  Special Functions and Their Applications. Englewood Cliffs, N.J. : Prentice-Hall, INC; 1965.

\end{enumerate}

\end{document}